\documentclass[12pt]{article}
\usepackage{amsfonts,latexsym,amssymb}
 \newcommand{\beqn}{\begin{eqnarray}}
 \newcommand{\eeqn}{\end{eqnarray}}
 \newcommand{\be}{\begin{equation}}
 \newcommand{\ee}{\end{equation}}
 \newcommand{\ba}{\begin{array}}
 \newcommand{\ea}{\end{array}}
 \newcommand{\pa}{\partial}
 \newcommand{\re}{\ref}
 \newcommand{\ci}{\cite}
 \newcommand{\la}{\label}

\newcommand{\Om}{\Omega}

\newcommand{\na}{\nabla}
\newcommand{\om}{\omega}
\newcommand{\bt}{\beta}
\newcommand{\al}{\alpha}

 \newcommand{\De}{\Delta}

\def\R{{\rm I\kern-.1567em R}}
\def\M{{\rm I\kern-.1567em M}}

\def\div {{\rm div}}

\def\dist{{\rm dist}}

 \newtheorem{theorem}{Theorem}[section]

\begin{document}

\begin{center} {\Large\bf Navier-Stokes equations: almost
$L_{3,\infty}$-case}\\
  \vspace{1cm}
 {\large
  G. Seregin}

 \end{center}

 \vspace{1cm}
 \noindent
 {\bf Abstract } A sufficient condition of regularity for
 solutions to the Navier-Stokes equations is proved. It
 generalizes the so-called $L_{3,\infty}$-case.

 \vspace {1cm}

\noindent {\bf 1991 Mathematical subject classification (Amer.
Math. Soc.)}: 35K, 76D.


\setcounter{equation}{0}
\section{Introduction  }

The paper is motivated by the following result in \ci{S2}. Let us
consider two functions $v$ and $p$ defined in $Q_T=\Om\times
]0,T[$, where $\Om\subset\mathbb R^3$ and $T$ is a positive
parameter. Assume that they meet three conditions:

\be\la{i1}v\in L_{2,\infty}(Q_T)\cap W^{1,0}_2(Q_T),\quad p\in
L_\frac 32(Q_T);\ee

$v$ {\it and} $p$ {\it satisfy the Navier-Stokes equations}
\be\la{i2}\left.
\begin{array}{c}\pa_t v+v\cdot\na\,v -\Delta \,v =-\nabla
\,p  \\

            \div \,v = 0\end{array}
            \right\}\qquad\mbox{in}\,\,\, Q_T\ee

{\it in the sense of distributions};

$v$ {\it and} $p$ {\it satisfy the local energy inequality}
\begin{eqnarray}\la{i3}
 \nonumber  \int\limits_\Om\varphi(x,t)|v(x,t)|^2\,dx+2\int
 \limits_{\Om\times
]0,t[}\varphi |\na
v|^2\,dxdt^\prime &&  \\
\leq\int\limits_{\Om\times ]0,t[}(|v|^2(\De
\varphi+\pa_t\varphi)+v\cdot\na\varphi ( |v|^2+2p))\,dxdt^\prime
&&
\end{eqnarray}

{\it for a.a.}  $t\in ]0,T[$ {\it and for all nonnegative smooth
functions, vanishing

in the neighborhood of the parabolic boundary}  $\pa^\prime Q_T$
{\it of the space-time cylin-

der} $Q_T$.
%\equiv\Om\times
%\{t=0\}\cup\pa\Om\times [0,T]$.

A pair $v$ and $p$, having properties (\ref{i1})-(\ref{i3}), is
also called a suitable weak solution to the Navier-Stokes
equations in $Q_T$. Such kind of solutions was defined and treated
by Scheffer \ci{Sc}, \ci{Sc1}, Caffarelli-Kohn-Nirenberg \ci
{CKN}, and others (see, for instances, \ci{Li}, \ci{LS}, \ci{S3},
and \ci{S4}). Our version of the definition of suitable weak
solutions is due to F.-H. Lin {\ci{Li}. It seems to be more
convenient to study.

We also say that the space-time point $z_0=(x_0,t_0)$ with $x_0\in
\Om$ and $0<t_0\leq T$ is a regular point of $v$ (in the
Caffarelli-Kohn-Nirenberg sense) if there exists a positive number
$r<\min\{\dist\{x_0,\Om\}, \sqrt{t_0}\}$ such that $v\in
L_\infty(Q(z_0,r))$, where $Q(z_0,r)=B(x_0,r)\times ]t_0-r^2,t_0[$
and $B(x_0,r)$ is the three-dimensional ball
%in $\mathbb R^3$
of
radius $r$ with the center at the point $x_0$. The point $z_0$ is
called singular if it is not regular.

We may give different definitions of regular points. For example, we
can replace the space $L_\infty(Q(z_0,r))$ with
$C(\overline{Q}(z_0,r))$ (or even with $C^\al(\overline{Q}(z_0,r))$
for some positive $\al$). However, according to the local regularity
theory for the Stokes equations, they turned out to be equivalent.

We may ask the following question: how many singular points are,
say, at $t=T$? In \ci{S2}, it was shown that there exists a
positive  universal constant $\varepsilon$ such that \be\la{i4}
N_T\leq \varepsilon\limsup\limits_{t\uparrow T}\frac 1{T-t}\int
\limits_t^T\int\limits_{\Om}|v(x,s)|^3\,dxds.\ee Here, $N_T$ is
the number of singular points at $t=T$. Obviously, if the right
hand side of (\ref{i4}) is finite, then  $N_T$ is finite.  For
example, this can happen if \be\la{i5}v\in L_{3,\infty}(Q_T).\ee
The similar result  was established earlier by Neustupa \ci{Ne}.
Later, in \ci{SS1}, \ci{ESS4}, it was proved that (\ref{i5})
implies regularity of $v$ in $Q_T$ and thus $N_T=0$. Now, it is
interesting to figure out what happens if \be\la{i6}
M_T\equiv\limsup\limits_{t\uparrow T}
%\sup\limits_{0<t< T}
\frac
1{T-t}\int
\limits_t^T\int\limits_{\Om}|v(x,s)|^3\,dxds<+\infty.\ee

Clearly, (\ref{i6}) is less restrictive than (\ref{i5}).
Nevertheless, the main result of the paper says that the answer is
the same.
\begin{theorem}\la{it1} Assume that $v$ and $p$ satisfy conditions
(\ref{i1})-(\ref{i3}).
%, (\ref{i6}).
Let, in addition, \be\la{i7} m_T =\liminf\limits_{t\uparrow
T}\frac 1{T-t}\int
\limits_t^T\int\limits_{\Om}|v(x,s)|^3\,dxds<+\infty.\ee Then,
points $z=(x,T)$, where $x\in\Om$, are regular.\end{theorem}

To demonstrate how  Theorem \ref{it1} can be used, let us consider
the following Cauchy problem for the Navier-Stokes equations

\be\la{a1}\left.
\begin{array}{c}\pa_t v+v\cdot\na\,v -\Delta \,v =-\nabla
\,p  \\

            \div \,v = 0\end{array}
            \right\}\qquad\mbox{in}\,\,\, Q_\infty,\ee
where $Q_\infty=\mathbb R^3 \times ]0,+\infty[$,
%$\Om$  is a
%bounded domain in $\mathbb R^3$, \be\la{a2} v(x,t)=0, \qquad
%x\in\pa\Om,\quad 0\leq t< \infty,\ee
\be\la{a3}v(x,0)=a(x),\qquad x\in \Om.\ee It is assumed that
\be\la{a4}a\in C^\infty_0(\mathbb R^3), \qquad \div\,a=0.\ee

It is well known that problem (\ref{a1})-(\ref{a3}) has at least one
weak solution which is called the  weak Leray-Hopf solution (see
%,
%for example,
monograph \ci{L} for details). One of the challenging problem in
mathematical hydrodynamics is to show that the above solution is
unique. In turn, among  various approaches to this problem, the
idea to prove smoothness of the weak Leray-Hopf solutions is quite
popular. The following theorem might be regarded as a small step
in that direction.
\begin{theorem}\la{at1} Let us denote by $T$ the first moment of time
when singular points appear. %If $T<+\infty$, t
Then
$$\lim\limits_{t\uparrow T}\frac 1{T-t}\int\limits^T_t\int
\limits_{\mathbb R^3}
|v|^3dz=+\infty.$$
\end{theorem}
\textsc{Proof} Theorem \re{at1} is a direct consequence of Theorem
\ref{it1}.

\setcounter{equation}{0}
\section{Proof of Theorem \ref{it1} }

By the natural scaling for the Navier-Stokes equations, it is
sufficient to replace the cylinder $Q_T$ with the cylinder
$Q=B\times ]-1,0[$ and to prove that $z=0$ is regular point. Here,
$B(r)=B(0,r)$ and $B=B(1)$. Now, condition (\ref{i7})
%takes
may be taken in the form \be\la{p1} M\equiv\liminf\limits_{t\uparrow
0 }\,-\frac 2{t}\int
\limits_t^0\int\limits_{B}|v(x,t)|^3\,dxds<+\infty.\ee

Next, by the known multiplicative inequality \be\la{p2} v\in
L_{\frac {10}3}(Q),\ee\be\la{p3}|v|\,|\na\,v|\in L_{\frac 98,\frac
32}(Q),\ee and since \be\la{p4}p\in L_\frac 32(Q),\ee one can
apply the local regularity theory for the Stokes system and
conclude that \be\la{p5}|\pa_t v|+|\na^2v|+|\na p|\in L_{\frac
98,\frac 32}(Q(5/6)).\ee Here, $Q(r)=Q(0,r)$ and $Q=Q(1)$. In
particular, (\ref{p5}) allows us to fix a representative of the
function $t\mapsto v(\cdot,t)$ in such a way that $$the \,\,
function\,\,t \mapsto \int\limits_{B(5/6)}v(x,t)\cdot w(x)\,dx\,\,
is\,\, continuous\,\, on \,\, [-(5/6)^2,0]\,\, $$$$for \,\,any
\,\,w\in L_2(B(5/6))$$ and thus, for each $t\in [-(5/6)^2,0]$,
\be\la{p6}\|v(\cdot,t)\|_{L_2(B(5/6))}<+\infty.\ee

Now, our aim is to show that \be\la{p7}v(\cdot,0)\in
L_3(B(5/6)).\ee To this end, we note that, by (\ref{p1}), there
exists a sequence $t_k\in ]-1,0[$ such that $t_k\uparrow 0$ as
$k\to \infty$ and, for $k=1,2,...$, \be\la{pd1}-\frac
1{t_k}\int\limits^0_{t_k}\int\limits_B|v(x,t)|^3dxdt\leq M.\ee
Then, we introduce the additional notation
$$g(t)=\|v(\cdot,t)\|_{L_3(B(5/6))},$$$$E_k=\{t_k<t<0\,\,\|\,\,g(t)>10M\},\qquad
E'_k=]t_k,0[\setminus E_k,$$ where $k=1,2,...$. By (\ref{pd1}), we
find
$$-\frac
1{t_k}|E_k|10M
%\leq -\frac
%1{t_k}\int\limits^{t_{k+1}}_{t_k}g(t)\,dt
\leq -\frac 1{t_k}\int\limits^{0}_{t_k}g(t)\,dt\leq M.$$ So,
$$|E_k|\leq \frac {|t_k|}{10}$$ and thus $$|E'_k|=|t_k|-|E_k|
\geq\frac 9{10}|t_k|>0.$$ Therefore, for each $k=1,2,...$, there
exists $s_k\in ]t_k,0[$ such that \be\la {p8}s_k\uparrow 0 \,\, as
\,\, k\to \infty\qquad and\qquad g(s_k)\leq 10M.\ee On the other
hand, according partial regularity theory of the Navier-Stokes
equations
$$v(x,s_k)\to v(x,0), \qquad \forall x \in B\setminus\Sigma,$$ and
the 1D Hausdorff measure of $\Sigma$ is zero. By Fatou's lemma and
by (\ref{p8}), we have
$$\int\limits_{B(5/6)}|v(x,0)|^3dx\leq 10M.$$
So, (\ref{p7}) is proved.

Now, we split the pressure into two parts: \be\la{p9}p=p^1+p^2,\ee
where $p^1$ is determined as a unique solution of the following
problem $$\int\limits_Bp^1(x,t)\Delta
\varphi(x)\,dx=-\int\limits_Bv(x,t)\otimes v(x,t):\na^2
\varphi(x)\,dx.$$ Here, $\varphi$ is an arbitrary  test function
from the space $\{\varphi\in W^2_3(B)\,\,\|\,\,\varphi=0\,\,on
\,\,\pa B\}$. It is well known that $p^1$  satisfies the estimate
\be\la{p10}\int\limits_B|p^1(x,t)|^\frac 32\,dx\leq c
\int\limits_B |v(x,t)|^3\,dx\ee and thus \be\la{p11} -\frac
1{t_k}\int\limits^0_{t_k}\int\limits_B|p^1(x,s)|^\frac
32\,dxds\leq cM,\qquad k=1,2,....\ee The second component of the
pressure is a harmonic function and, therefore, satisfies the
estimates:
$$\sup\limits_{x\in B(2/3)}|p^2(x,t)|^\frac 32\leq
c\int\limits_B|p^2(x,t)|^\frac 32\,dx$$ \be\la{p12}\leq
c\Big(\int\limits_B|p(x,t)|^\frac
32dx+\int\limits_B|p^1(x,t)|^\frac 32dx\Big)\ee
$$\leq c\Big(\int\limits_B|p(x,t)|^\frac
32dx+\int\limits_B|v(x,t)|^ 3dx\Big).$$

Now, we extend functions $v$, $p^1$, and $p^2$ by zero to the
whole space $\mathbb R^3\times \mathbb R^1$.

Assume that the statement of the theorem is false. Then, we know
that
%there exists a sequence $\{R_k\}^\infty_{k=1}$ such that
%$R_k\downarrow 0$ as $k\to +\infty$ and
 \be\la{p13} \frac
1{R^2}\int\limits_{Q(R)}\Big(|v|^3+|p|^\frac 32\Big)\,dz\geq
\varepsilon>0, \qquad\forall R\in ]0,1]\ee for some positive
universal constant $\varepsilon$. Fixing $T<-1000$, we may blow up
our solution at zero with the help of  the following scaling:
$$u^k(y,s)=R_kv(R_ky,R^2_ks),  \qquad (y,s)\in \mathbb R^3\times
\mathbb R^1,$$
$$q^{1k}(y,s)=R^2_kp^1(R_ky,R^2_ks),\qquad
q^{2k}(y,s)=R^2_kp^2(R_ky,R^2_ks),$$ where $R_k=\sqrt{t_k/T}\to 0$
as $k\to \infty $. We have (remember $t_k=TR^2_k$)
$$-\frac 1T\int\limits^0_T\int\limits_{\mathbb
R^3}|u^k(y,s)|^3dyds=-\frac
1T\int\limits^0_T\int\limits_{B(1/R_k)} |u^k(y,s)|^3dyds$$
\be\la{p14}= -\frac
1{TR^2_k}\int\limits^0_{TR^2_k}\int\limits_{B}|v(y,s)|^3dyds\leq
M.\ee
%provided \be\la{p15}-TR^2_k<1.\ee
$q^{1k}$ can be treated in
the similar way:
$$-\frac
1T\int\limits^0_T\int\limits_{\mathbb R^3}|q^{1k}(y,s)|^\frac
32dyds=-\frac 1T\int\limits^0_T\int\limits_{B(1/R_k)}
|q^{1k}(y,s)|^\frac 32dyds$$ \be\la{p16}= -\frac
1{TR^2_k}\int\limits^0_{TR^2_k}\int\limits_{B}|p^1(y,s)|^\frac
32dyds\leq cM.\ee
%provided that (\ref{p15}) holds.

As to $q^{2k}$, we take into account (\ref{p12}) and argue as
follows:
$$-\frac 1T\int\limits^0_T\int\limits_{B(a)}|q^{2k}(y,s)|^\frac
32dyds=-\frac
1{TR^2_k}\int\limits^0_{TR^2_k}\int\limits_{B(aR_k)}|p^2(y,s)|^\frac
32dyds$$
$$\leq -c\frac {(aR_k)^3}{TR^2_k}\int\limits^0_{TR^2_k}
\int\limits_{B}\Big (|p(x,t)|^\frac 32+|v(x,t)|^3\Big)dxdt$$
\be\la {p17}\leq -c\frac {a^3}TR_k\int\limits_Q\Big (|p|^\frac
32+|v|^3\Big)dz\to 0\ee as $k\to +\infty$ if $aR_k<2/3$.

Selecting subsequences (still denoted in the same way), we
have$$u^k\rightharpoonup u \qquad in\quad L_3(B(a)\times ]T,0[),$$
\be\la{p18}q^{1k}\rightharpoonup q \qquad in\quad L_\frac
32(B(a)\times ]T,0[),\ee
$$q^{2k}\rightarrow 0\qquad in\quad L_\frac
32(B(a)\times ]T,0[)$$ %for any $T<0$ and
for any $a>0$. Moreover, \be\la{p19} |u|^2 +|q|\in L_\frac
32(\mathbb R^3\times]T,0[).\ee
%for any $T<0$.

It remains to show that the pair $u$ and $q$ is a suitable weak
solution to the Navier-Stokes equations on  sets of the form
$B(a)\times ]T,0[$. To this end, we first observe that
\be\la{p20}\|u^k\|_{L_{2,\infty}(B(a)\times ]T,0[)}+ \|\na
u^k\|_{L_{2}(B(a)\times ]T,0[)}\leq C(a,T,M)<+\infty\ee and,
therefore, by known multiplicative inequalities,
\be\la{p21}\||u^k|\,|\na u^k|\|_{L_{\frac 98,\frac 32}(B(a)\times
]T,0[)}+ \| u^k\|_{L_{\frac {10}3}(B(a)\times ]T,0[)}\leq
C(a,T,M).\ee Next, the linear theory says that \be\la{p22}\||\pa_t
u^k|+|\na^2 u^k|\|_{L_{\frac 98,\frac 32}(B(a)\times ]T,0[)}\leq
C(a,T,M).\ee Estimates (\ref{p20})-(\ref{p22}), together with
known compactness arguments, imply
$$u^k\rightarrow u\qquad in\quad L_3(B(a)\times ]T,0[),$$
\be\la{p23}u^k\rightarrow u\qquad in\quad C([T,0];L_\frac
98(B(a)))\ee
%for any $T<0$ and
for any $a>0$.

Now, we can pass to the limit in the Navier-Stokes equations and
in the local energy inequality for $u^k$ and $q^k$ on sets of the
form $B(a)\times ]T,0[$ and conclude that limit functions $u$ and
$q$ generate a suitable weak solution to the Navier-Stokes
equations on those sets. The function $u$ has the  properties
$$u\in L_{2,\infty}(B(a)\times ]T,0[)\cap W^{1,0}_2(B(a)\times ]T,0[)
\cap L_\frac {10}3(B(a)\times ]T,0[),$$ \be\la{p24}|u|\,|\na u|+
|\pa_tu|+|\na^2 u|\in L_{\frac 98,\frac 32}(B(a)\times ]T,0[),\ee
$$u\in C([T,0];L_\frac {9}8(B(a))).$$

Let us show that our blow-up solution is not trivial. By scaling
and by (\ref{p13}), \be\la{pd2}\frac
1{(aR_k)^2}\int\limits_{Q(aR_k)}\Big(|v|^3+|p|^\frac
32\Big)dz=\frac 1{a^2}\int\limits_{Q(a)}\Big(|u^k|^3+|q^k|^\frac
32\Big)dz\geq \varepsilon>0\ee for all $a\in ]0,1]$ and for all
$k=1,2,...$. It follows from (\ref{p23}) that \be\la{pd3}
\int\limits_{Q(a)}|u^k|^3dz\rightarrow
\int\limits_{Q(a)}|u|^3dz.\ee Going back to the definition of
$p^1$, we find after change of variables:
$$\int\limits_{B(1/R_k)}q^{1k}(y,s)\De \psi(y)dy=
-\int\limits_{B(1/R_k)}u^k(y,s)\otimes u^k(y,s):\na^2\psi(y)dy$$
for any test function $\psi\in W^2_3(B(1/R_k))$ with $\psi=0$ on
$\pa B(1/R_k)$ and for any $s\in [-1/R^2_k,0]$. Next,
%Then we split
$q^{1k}$ can be split in the following way:
\be\la{pd4}q^{1k}=r^{1k}+r^{2k}\qquad in \quad Q(2),\ee where
$$\int\limits_{B(2)}r^{1k}(y,s)\De \psi(y)dy=
-\int\limits_{B(2)}u^k(y,s)\otimes u^k(y,s):\na^2\psi(y)dy$$ for any
test function $\psi\in W^2_3(B(2))$ with $\psi=0$ on $\pa B(2)$ and
for any $s\in [-2^2,0]$. For $r^{1k}$, we have the estimate
\be\la{pd5}\int\limits_{Q(2)}|r^{1k}|^\frac 32 dz\leq
c\int\limits_{Q(2)}|u^k|^3dz \rightarrow
c\int\limits_{Q(2)}|u|^3dz.\ee From (\ref{pd4}), it follows that
$$\De r^{2k}(\cdot,s)=0\qquad in \quad B(2)$$ for any
$s\in [-2^2,0]$. By properties of harmonic functions,
$$\sup\limits_{y\in B(1)}|r^{2k}(y,s)|^\frac 32\leq
c\int\limits_{B(2)} |r^{2k}(y,s)|^\frac 32dy, \qquad s\in
[-2^2,0].$$ So, letting
$e=(y,s)$,$$\int\limits^0_{-2^2}\|r^{2k}(\cdot,s)\| ^\frac
32_{L_\infty(B(1))}ds\leq c\int\limits_{Q(2)} |r^{2k}|^\frac
32de\leq $$
$$\leq  c\int\limits_{Q(2)} |q^{1k}|^\frac 32de
+c\int\limits_{Q(2)} |r^{1k}|^\frac 32de$$ \be\la{pd6}\leq
c\int\limits_{Q(2)} |q^{1k}|^\frac 32de +c\int\limits_{Q(2)}
|u^{k}|^3de\leq C(M,T).\ee Then, we have (see (\ref{pd2}) and
(\ref{pd3}))
$$0<\varepsilon=\limsup\limits_{k\to\infty}\,\frac 1{a^2}\int
\limits_{Q(a)}\Big(|u^k|^3+|q^k|^\frac 32\Big)de$$
$$\leq \frac 1{a^2}\int
\limits_{Q(a)}|u|^3de+c\limsup\limits_{k\to\infty}\,\frac
1{a^2}\Big(\int \limits_{Q(a)}|q^{1k}|^\frac 32de+\int
\limits_{Q(a)}|q^{2k}|^\frac 32de\Big).$$ By the last relation in
(\ref{p18}), we have
$$0<\varepsilon\leq \frac 1{a^2}\int
\limits_{Q(a)}|u|^3de+c\limsup\limits_{k\to\infty}\,\frac 1{a^2}
\int \limits_{Q(a)}|q^{1k}|^\frac 32de.$$ Taking into account
(\ref{pd4})--(\ref{pd6}), we find
$$0<\varepsilon\leq \frac 1{a^2}\int
\limits_{Q(a)}|u|^3de+c\limsup\limits_{k\to\infty}\,\frac 1{a^2}
\int \limits_{Q(a)}|r^{1k}|^\frac 32de$$
$$+c\limsup\limits_{k\to\infty}\,\frac 1{a^2}
\int \limits_{Q(a)}|r^{2k}|^\frac 32de$$
$$\leq \frac 1{a^2}\int
\limits_{Q(a)}|u|^3de+ \frac c{a^2}\int \limits_{Q(2)}|u|^3de$$
$$+\frac 2{a^2}|B(a)|\int\limits^0_{-a^2}
\sup\limits_{y\in B(1)}|r^{2k}(y,s)|^\frac 32ds$$
$$\leq \frac c{a^2}\int \limits_{Q(2)}|u|^3de+C(M,T)a.$$
Choose $a$ sufficiently small so that $$C(M,T)a\leq \frac 12
\varepsilon.$$ Then \be\la{pd7}0<\varepsilon a^2\leq c\int
\limits_{Q(2)}|u|^3de.\ee So, $u$ is not trivial.

%Moreover, $u$ is not trivial. Indeed,  \be\la{p25} \frac
%1{R^2_k}\int\limits_{Q(R_k)}|v|^3\,dz=\int\limits_{Q}|u^k|^3\,dz
%\rightarrow\int\limits_{Q}|u|^3\,dz\geq \varepsilon>0.\ee

 Now, we are in a position to show that
\be\la{p26}u(\cdot,0)=0\qquad in \quad \mathbb R^3.\ee To this
end, we proceed as follows. For any $a>0$, we have
$$\frac 1{a^2}\int\limits_{B(a)}|u(x,0)|\,dx\leq\frac 1{a^2}
\int\limits_{B(a)}|u^k(x,0)-u(x,0)|\,dx$$$$ +\frac
1{a^2}\int\limits_{B(a)}|u^k(x,0)|\,dx=\bt_k+\frac
1{a^2}\int\limits_{B(a)}|u^k(x,0)|\,dx$$$$=\bt_k+\frac
1{(aR_k)^2}\int\limits_{B(aR_k)}|v(x,0)|\,dx$$
$$\leq\bt_k+c\Big(\int\limits_{B(aR_k)}|v(x,0)|^3\,dx\Big)^\frac
13.$$ By (\ref{p7}) and by (\ref{p23}), the right hand side of the
latter inequality tends to zero as $k\to+\infty$. So, (\ref{p26})
is proved.

Other important step in our scheme is to show that, given $T\in
]-\infty,0[$,  there exists $R(T)\in ]0,+\infty[$ such that
functions $u$ and $q$ are smooth in spatial variables  in
$(\mathbb R^3\setminus B(R(T)))\times [T/2,0]$. To this end, it is
sufficient to prove that \be\la
{p27}\int\limits_{Q(z_0,1)}(|u|^3+|q|^\frac 32)dz\rightarrow
0,\qquad z_0=(x_0,t_0),\ee as $|x_0|\to +\infty$ and $T/2\leq
t_0\leq 0$. In fact, (\ref{p27}) follows from (\ref{p19}) (see
similar arguments in \ci{ESS4}). Then, according to Lemma 2.2 in
\ci{ESS4}, we can state  that \be\la{p28}|u|+|\na u|\leq
C<+\infty\qquad in \quad (\mathbb R^3\setminus B(R(T)))\times
[T/2,0].\ee As it was shown in \ci{ESS4}, (\ref{p26}) and
(\ref{p28}) implies that
$$\om(x,t)=0, \qquad x\in \mathbb R^3,\quad T/2\leq t\leq 0,$$
where $\om=\na \wedge u$ is the vorticity. So, $u(\cdot,t)$ is a
harmonic function in $\mathbb R^3$ for all $T/2\leq t\leq 0$. On
the other hand, it follows from (\ref{p19}) that, for a. a.
$T/2\leq t\leq 0$, $u(\cdot,t)$ is in $L_3(\mathbb R^3)$ and,
therefore, $u(\cdot,t)=0$ for the same $t$. But this is in a
contradictions with (\ref{pd7}). Theorem \ref{it1} is proved.

\noindent

G. Seregin\\
Steklov Institute of Mathematics at St.Petersburg, \\
St.Peterburg, Russia

\end{document}